\documentclass[a4paper,11pt]{amsart}
\usepackage{amssymb}
\usepackage{amscd}

\def\Box{\square}
\def\ci{C^\infty}
\def\circe{{\scriptstyle{\circ}}}
\def\cph{{\mathbb{CP}}^2}
\def\cM{{\mathcal{M}}}

\def\Bbb{\mathbb}
\def\Cal{\mathcal}

\def\IQ{{\Cal{I}}^{{\rm Q}}}

\def\bb{{\bf b}}
\def\bbC{{\bf C}}
\def\bbP{{\bf P}}

\let\b=\beta
\let\k=\kappa
\let\l=\lambda

\def\kth{k^{\underline{{\rm th}}}}
\def\lth{\ell^{\underline{{\rm th}}}}
\def\slp{{\rm(}}
\def\srp{{\,\rm)}}

\def\bbR{{\mathbb{R}}}
\def\bbT{{\mathbb{T}}}

\def\dint{{\displaystyle\int}}

\let\d=\delta
\let\e=\varepsilon

\let\i=\iota

\let\w=\omega
\let\r=\rho

\def\LOT{{\rm LOT}}

\newcommand{\IT}[1]{{\rm(}{\it{\!#1}}{\rm)}}

\let\ccdot\cdot
\def\cdot{\hbox to 2.5pt{\hss$\ccdot$\hss}}

\newcommand{\om}{\omega}
\newcommand{\si}{\sigma}
\newcommand{\Om}{\Omega}

\newtheorem{theorem}{Theorem}[section]
\newtheorem{lemma}[theorem]{Lemma}
\newtheorem{proposition}[theorem]{Proposition}
\newtheorem{corollary}[theorem]{Corollary}

\newcommand{\cC}{{\Cal C}}
\newcommand{\cH}{{\Cal H}}
\newcommand{\cN}{{\Cal N}}
\newcommand{\cV}{{\Cal V}}
\newcommand{\ce}{{\Cal E}}
\newcommand{\cE}{{\Cal E}}
\newcommand{\cQ}{{\Cal Q}}
\newcommand{\cK}{{\Cal K}}
\newcommand{\cR}{{\Cal R}}
\newcommand{\nd}{\nabla}
\newcommand{\Rho}{{\mbox{\sf P}}}
\newcommand{\End}{\operatorname{End}}
\def\cI{{\Cal I}}
\def\cIQ{\cI^{{\rm Q}}}
\newcommand{\wh}{\widehat}

\newcommand{\nn}[1]{(\ref{#1})}
\newcommand{\nnn}[1]{{\rm(}\!{\ref{#1}}{\rm)}}
\newcommand{\bT}{\mbox{{$\Bbb T$}}}
\newcommand{\bg}{\mbox{\boldmath{$ g$}}}

\newcommand{\V}{{\mbox{\sf P}}}
\newcommand{\J}{{\mbox{\sf J}}}

\begin{document}
\renewcommand{\today}{}
\title{Pontrjagin forms and invariant objects related to the Q-curvature}
\author{Thomas Branson and A. Rod Gover}

\begin{abstract}
We clarify the conformal invariance of the Pontrjagin forms by giving
them a manifestly conformally invariant construction; they are shown
to be the Pontrjagin forms of the conformally invariant tractor
connection. The Q-curvature is intimately related to the
Pfaffian. Working on even-dimensional manifolds, we show how the
$k$-form operators $Q_k$ of \cite{ddd}, which generalise the
Q-curvature, retain a key aspect of the $Q$-curvature's relation to
the Pfaffian, by obstructing certain representations of natural
operators on closed forms. In a closely related direction, we show
that the $Q_k$ give rise to conformally invariant quadratic forms
$\Theta_k$ on cohomology that interpolate, in a suitable sense,
between the integrated metric pairing (at $k=n/2$) and the Pfaffian
(at $k=0$). 
Using a different construction, we show that the $Q_k$
operators yield a generalisation of the period map which maps 
conformal structures to Lagrangian subspaces of the direct sum
$H^k\oplus H_k$ (where $H_k$ is the dual of the de Rham cohomology
space $H^k$).  We couple the $Q_k$ operators with the Pontrjagin forms
to construct new natural densities that have many properties in common
with the original Q-curvature; in particular these integrate to global
conformal invariants.  We also work out a relevant example, and show
that the proof of the invariance of the (nonlinear) action functional
whose critical metrics have constant Q-curvature extends to the action
functionals for these new Q-like objects. Finally we set up eigenvalue
problems that generalise to $Q_k$-operators the Q-curvature
prescription problem.
\end{abstract}

\maketitle \markboth{{\sf Branson \& Gover}}{{\sf Pontrjagin
forms, Q-space}}

\section{Introduction}
\subsection{Overview}
In recent years there has been much interest in Q-curvature on
even-dimensional conformal manifolds
\cite{BCY,CGY,CQY,CYAnnals,CY,FGrQ,FHmrl,GoPetLap,RobinKengo,GrZ,gursk,gurskyv}.
Q-curvature naturally appears in {\em Polyakov-type formulas} for the
determinant or torsion quotient with respect to a differential
operator or elliptic complex with good conformal behavior.  It is
shown in \cite{RobinKengo} that the integral of the Q-curvature is
the action functional for the {\em Fefferman-Graham obstruction
tensor}.  This is in turn associated with the AdS/CFT
correspondence, and in fact 
quantities having some properties in common with the 
Q-curvature appear in the 
volume asymptotics of conformally compact Poincare-Einstein spaces.  
In terms of elementary geometric analysis, the prescription problem
for the Q-curvature provides a higher-dimensional generalisation
of the 2-dimensional Gauss curvature prescription problem.

The above-mentioned Polyakov-type formulas have motivated study of the
structure of natural densities $U$ having conformally invariant
integral.  {\em Conformal indices} in the sense of \cite{bocm}, or
{\em trace anomalies} in Physics, are quantities built from spectral
zeta functions $\zeta(s)$ (and specifically from their conformally
invariant values $\zeta(0)$ at $s=0$); these take the form $\int U$.
The {\em functional determinant} or {\em one-loop effective action} is
a spectral invariant that provides a regularisation of a divergent
{\em functional integral}; a torsion quantity is a certain well-chosen
linear combination of functional log-determinants $\zeta'(0)$.  The
local part of the conformal variation of these determinants are given
by integrals of quantities $U$ as above against arbitrary functions
(conformal factors).  This makes it worthwhile to understand the space
of natural densities $U$ that can appear; see e.g.\ \cite{ds}.  Though
this is easy to do in, for example, dimension 4, by writing out all
the possibilities, the situation is not well understood in higher
dimensions, despite the relevance of this problem to string and brane
theories; see \cite{erd} for some work in dimension 8. A conjecture of
central importance is the proposal that each natural
density $U$ with conformally invariant integral is the sum of a
multiple of the Pfaffian, a natural exact divergence, and a natural
locally conformally invariant density. Settling this is proving to be
not at all straightforward, see \cite{Spiros} for some recent
progress. 
 The Q-curvature captures the conformal content carried by
the Pfaffian, in the sense that its integral is a fixed multiple of
the Euler characteristic in the locally conformally flat case, so this
conjecture is equivalent to the analogous one in which the Pfaffian is
replaced by the Q-curvature.  A separate conjecture \cite{tor} asserts
that the Q-curvature and suitable modifications are the only terms
that can appear in a certain way in the local part of a Polyakov-type
formula.

A central result here is a new construction of natural densities that
are ``Q-like'' in the sense that each of these densities integrates
(in the compact setting) to a conformal invariant and has conformal
variation by a formally self-adjoint linear operator on the log of the
conformal factor. There has been recent progress in the construction
of such densities \cite{FHmrl,GoPetLap}. The construction here in
Theorem \ref{naturalQs} takes a completely new direction in that the
components of the construction are essentially non-local; they arise
from coupling Pontrjagin forms with the local expression for a global
(integrated) conformal pairing.

Our analysis of the structure of the de\,Rham complex from the
viewpoint of conformal structure in \cite{ddd} has uncovered a series
of linear differential operators $Q_k$ on lower degree cocyclic forms
that admit linear conformal change laws -- that is, which ``act like
Q-curvatures''.  In fact, the classical Q-curvature $Q$ of
\cite{tbkorea} is essentially the lowest-valence member of this series
($Q=Q_0 1$).  The $Q_k$ compress to maps between cohomologies and
spaces of conformal harmonics and this generalises the fact that $\int
Q_0$ is a numerical conformal invariant.  One of the key results is
Theorem \ref{nogo} which proves that these $Q_k$-operators encode
global information.  The conformal change of $Q$ is by an exact divergence, yet
$Q$ is not just the sum of a local conformal invariant and an exact
divergence. In the conformally flat case at least, the Pfaffian is the
obstruction. Theorem \ref{nogo} shows the analogue is true for the
$Q_k$ operators. The $Q_k$ changes conformally by an operator with a
left composition factor of the differential form coderivative $\d$,
yet $Q_k$ is not just the sum of a conformally invariant operator and
something with a left $\d$ factor.  This establishes that Q-curvatures
and Q-operators encode deep information about conformal structures
beyond what is available from the classification of locally
conformally invariant tensors and operators.  In particular, while
they do give conformally invariant cohomology maps these do not arise
from objects with local conformal invariance.  This is an indication
that the $k$-form Q operators for $k\ne 0$ are {\em conformal}
analogues of the (scalar) Pfaffian.

The global conformal
content of the $Q_k$-operators is captured
by a global pairing in Theorem \ref{capform}. This shows that in a
suitable sense and under suitable conditions on the conformal
structure (conditions which we expect to be generic), the $Q_k$ give
rise to conformally invariant real-valued quadratic forms $\Theta_k$
on the cohomologies $H_k$. When $k=n/2$, this quadratic form is the
integrated metric pairing, while for connected, conformally flat
manifolds, $\Theta_0$ is a multiple of the Euler characteristic.
(Relaxing the connectedness assumption, it multiplies the class in
0-cohomology contributed by a given connected component by that
component's Euler characteristic.)  Thus one may view the $Q_k$
operators as special local and natural expressions for these global
quantities. 
In Section \ref{moments} we establish eigenvalue problems
that generalise to $Q_k$-operators the
Q-curvature prescription problem.  For a class of 10-dimensional manifolds,
we exhibit a solution to such a prescription problem 
in Section \ref{cpH}.  These examples establish that 
the new objects constructed in Theorem \ref{naturalQs} are not trivial
in general.

It should be interesting to determine whether the $\Theta_k$ distinguish 
conformally flat structures on a given
fixed manifold or whether, in this setting, they reflect only
topological data.  On the other hand in general we expect the $Q_k$
operators to be sensitive to conformal structure and in Proposition
\ref{periodprop} we construct a map from the space of conformal structures on
$M$ to the Grassmannian of Lagrangian subspaces of $H^k\oplus
H_k$. This generalises a special case of this generalised period map
developed by Eastwood and Singer \cite{EastSin2}.
For some applications of period maps to gauge theory, see
\cite{DonKron}, Section 4.3, and \cite{LeBrun,IL}.

Apart from those mentioned already, the main results are as follows.
We show that the Pontrjagin classes of the standard conformal
tractor bundle agree with the Pontrjagin classes for $TM$.  This
leads to a new and simple proof that the Pontrjagin classes of the
tangent bundle obstruct conformally flat metrics. The argument
uses no explicit information about forms representing the classes.
We then show that in fact the Pontrjagin {\em forms} for the
tractor bundle (not just their cohomology classes) agree with the
usual Pontrjagin forms. 
This provides a transparent explanation for the conformal invariance 
of these forms (first proved by Chern and Simons \cite{CS74}). It 
also suggests the point of view that the Pontrjagin forms are more
naturally associated to the tractor connection than to the
Levi-Civita connection, since the latter involves choices (in
particular a metric from the conformal class) that are not
necessary in order to write down these forms.

Theorem
\ref{couplethm} shows that the Q-operator of lowest positive (in
fact second) order is capable of producing Q-like objects even
when coupled to an auxiliary bundle and connection, while Theorem
\ref{nocouplethm} shows that the higher-order Q operators
definitely do not perform this service. Theorem \ref{invtfcnl}
establishes the conformal transformation invariance of the action
functional associated to the problem of prescribing constancy
for Q-like invariants.

We would like to thank Peter Gilkey and 
Claude LeBrun for helpful discussions,
and the National Center for Theoretical Sciences (Mathematics Division)
of Taiwan for support during a workshop where much of this work was done.
The second author would like to thank
the Royal Society of
New Zealand for support via Marsden Grant no.\ 02-UOA-108,
and the New Zealand Institute of Mathematics and its Applications
for support via a Maclaurin Fellowship.

\subsection{Basic constructions}
We define a vector space of Q-like objects, the {\em linear
Q-space}, as follows. Let $\cI$ be the space of natural
$(-n)$-densities $A$ on pseudo-Riemannian manifolds of even
dimension $n$; that is, quantities built in a universal and
polynomial way from a pseudo-Riemannian metric and its inverse, 
together with the 
associated covariant derivative $\nd$ and Riemann
curvature $R$.
Our normalisation of conformal density weights is such that
$(-n)$-densities are integrable given only a conformal structure;
henceforth we shall simply use the term {\em densities} to refer to
these.  In a natural density, for each monomial summand, the number of
$\nd$ plus twice the number of $R$ occurring is $n$.  Under a
conformal change of metric $\wh{g}=e^{2\w}g$, where $\w$ is a smooth
function, we have
$$
\wh{A}=A+A^1(d\w,g,g^{-1},\nd,R)+\cdots+A^n(d\w,g,g^{-1},\nd,R),
$$
where $A^s$ is a universal polynomial which is $s$-homogeneous in $d\w$.
It is an elementary exercise to show
that if $A^s$ is identically zero
(for all $g$ and $\w$), then so are the $A^t$ for
$t>s$. Thus $\cI^s:=\{A\mid A^{s+1}=0\}$ is a filtration of the
set of natural densities.  In particular, $\cI^0$ is the space of
conformally invariant natural densities.  It is also elementary to
show that if $A\in\cI^s$, then
\begin{equation}\label{invthg}
\w\mapsto A^s(d\w,g,g^{-1},\nd,R)
\end{equation}
is a conformally
invariant $s$-homogeneous differential operator from functions to
densities. In particular, if $s=1$, this is an invariant linear
differential operator.

The linear Q-space is a certain subspace of $\cI^1$, the densities
with a linear conformal change law.  The additional condition that
we shall impose reflects the origin of the Q-curvature concept in
the study of determinant and torsion quotient formulas
\cite{bo91}. Because the conformal variation of a Q-curvature
should be the second conformal variation of a (possibly nonlocal) quantity,
there should be a symmetry condition for this second variation.
This works out to be as follows.  For any $A\in\cI$, let $({\bf
b}A)\w:=A^1(d\w,g,g^{-1},\nd,R)$. Since ${\bf b}A$ is a linear
differential operator from functions to densities, and these are
naturally dual given a conformal structure, the formal adjoint
$(\bb A)^*$ is also a linear differential operator from functions
to densities. We define $\cI^{{\rm FSA}}$ to be the subspace of
$A\in\cI$ for which ${\bf b}A$ is formally self-adjoint. Note that
$d$ is a right composition factor for any $\bb A$.  Thus if
$A\in\cI^{{\rm FSA}}$, then (from the normal form of \cite{EGTN}) the real coefficient operator $\bb A$ has a
factorisation $\bb A=\d Nd$ for some natural operator $N$.

The linear Q-space is
$$
\cI^{{\rm Q}}:=\cI^1\cap\cI^{{\rm FSA}}.
$$
By the remark around \nn{invthg}, for each $Q\in\IQ$, the operator
${\bf b}Q$ is conformally invariant; we also know it is formally
self-adjoint
and of the form $\d Nd$.
In particular,
$$
\int\wh{Q}=\int(Q+\d Nd\w)=\int Q,
$$
since $\d Nd\w$ is an exact divergence.  This shows that $\int Q$
is conformally invariant.  Elements of $\IQ$ representing nonzero
classes modulo $\cI^0$ (i.e.\ elements which are not just local
conformal invariants) are especially interesting.

The linear Q-space makes precise a notion of ``Q-like''.
In this announcement, we indicate how one may allow the $Q_k$ to
act on characteristic forms to obtain elements of the linear Q-space.
Alternatively, instead of maximally reducing the valence in this way,
we may use characteristic forms to obtain additional operators on
cocyclic forms with linear conformal change laws.

\section{Conformal forms}

For a complex vector bundle $\cV$ of rank $N$, let us write $AB$ to
mean the composition of sections $A,B\in \Gamma(\End(V))$. By
construction the function
$$
s_k(A) = {\rm Tr}(\underbrace{AA\cdots A}_{k})
$$
is an invariant which is pointwise polynomial (and homogeneous of
degree $k$) in the $N^2$ components of $\End(V)$.
In the following, we shall call such objects {\em polynomial invariants}.
We recall that the
$s_k(A)$ polynomially generate all such invariants.

Since exterior multiplication is commutative on even forms, for a
given polynomial invariant $P(A)$ we can replace the section $A\in
\End(\cV)$ with a 2-form $F$ taking values in $\End(\cV)$ to obtain an
invariant, and in general non-trivial, $2k$-form that we denote
$P(F)$. That is, $P$ determines a map
$$
P:
\ce^2(\End(\cV))\to \ce^{2k},
$$
which is algebraic in the sections $F\in \ce^2(\End(\cV))$ since it
arises from a bundle map. This construction is most famliar when
$F=F^D$ is the curvature of some connection $D$ on $\cV$. In this case
the following results are well known (see e.g. \cite{KobNom}).

\begin{lemma}\label{charfacts}
For each polynomial invariant $P$ and connection $D$, $P(F^D)$ is a
 closed form. The de\,Rham cohomology class $[P(F^D)]\in H^*(M;{\Bbb
 C})$ is independent of the choice of connection $D$ on $\cV$, and
 isomorphic vector bundles define identical cohomology classes.
\end{lemma}

Recall that the $\lth$ Chern character class of the complex
vector bundle $\cV$ is the class in $H^{2\ell}(M;{\Bbb C})$ given by
\begin{equation}\label{cherndef}
c_\ell(\cV)=\frac{1}{\ell !}\left[s_\ell\left(
\frac{-1}{2\pi \sqrt{-1}} F^D\right)\right] ,
\end{equation}
for any connection $D$ on $\cV$. (It is slightly more convenient here
to use these character classes rather than the {\em Chern classes},
whose definition involves the elementary symmetric polynomial
$\sigma_\ell$ instead of $s_\ell$.)
The $\kth$ Pontrjagin
character class of a real vector bundle $\cV$ is simply
$c_{2k}(\cV_{\Bbb C})$, and takes values in real cohomology. In this
setting we will call the forms
\begin{equation}\label{pontdef}
p_k^D=\frac{1}{(2k)!}s_{2k}
\left(\frac{-1}{2\pi
\sqrt{-1}} F^D\right)
\end{equation}
{\em Pontrjagin forms}.

It is result of Chern and Simons \cite{CS74} that the Pontrjagin
forms for a Riemannian connection are conformally invariant. We
will show here that in fact these forms arise from the conformally
invariant tractor connection $\nabla$ \cite{BEGo,CapGoLuminy} on
the (real) tractor bundle $\Bbb T$. Let us write $\Omega$ for the
curvature of the standard tractor connection.  Since this
connection preserves a metric on $\bT$, the invariant $p$-forms
$P(\Omega)$ vanish identically unless they have degree $p=4k$.
(Contraction with several copies of the metric produces a self- or
skew-adjoint endomorphism according to whether $p$ is of the form
$4k$ or $4k+2$.)
\vspace{0.2cm}

\noindent{\bf Definition:} For each $k\geq 1$ we will call the closed $4k$-form
$
\tau_k = p_k^{\Bbb T}
$ the $\kth$ {\em conformal Pontrjagin \IT{character} form} for a
Riemannian or pseudo-Riemannian manifold and we will use the term {\em
conformal form} for any invariant form generated by the conformal
Pontrjagin (character) forms. We first observe that these recover the
usual character classes.
\vspace{0.2cm}

\begin{theorem}\label{tracclasses}
  For $k\geq 1$ the class $[\tau_k]\in H^{4k}(M)$ is ${\rm ph}_k(TM)$,
  the $k^{\rm th}$ Pontrjagin character class of the tangent bundle.
\end{theorem}

\noindent{\bf Proof:} From its definition \cite{BEGo,CapGoLuminy}
the tractor bundle $\bT$ is isomorphic (as a vector bundle,
without considering further structure) to the direct sum bundle
${\Bbb R}\oplus TM \oplus {\Bbb R}$. The latter may be equipped
with a connection which is a trivial extension of any connection
on $TM$, so the result follows from the definition and Lemma
\ref{charfacts}.$\qquad\Box$
\vspace{0.2cm}

Since the tractor curvature vanishes on manifolds that are
(locally) conformally flat, it follows that we have an alternative
proof of the following result of \cite{CS74}:

\begin{corollary}\label{obstructions}
The Pontrjagin character classes
$$
{\rm ph}_k(TM) \quad k\geq 1
$$
vanish on manifolds $M$ that admit a metric \slp of any signature\srp{}
which is
locally conformally flat.
\end{corollary}

Alternatively one obtains Theorem \ref{tracclasses} from the
stronger result that follows. Here $R$ denotes the Riemann
curvature tensor for Riemannian (or pseudo-Riemannian) structure
on $M$, and $C$ the corresponding Weyl curvature.

\begin{proposition}\label{confformeq}
On a manifold with a metric $g$ of any signature,
\begin{equation}\label{tCR}
\tau_k= \frac{1}{(2k)!}s_{2k}\left(\frac{-1}{2\pi \sqrt{-1}}
C\right)= \frac{1}{(2k)!}s_{2k}\left(\frac{-1}{2\pi \sqrt{-1}}
R\right),
\end{equation}
where in the middle expression we regard the Weyl curvature as an
$\End(TM)$-valued 2-form.
\end{proposition}

\noindent{\bf Proof:}
Following the notation of \cite{GoPetLap}, in terms of the metric
$g$ the tractor curvature is given by
$$
\Omega_{ab}{}^C{}_D
= Z^C{}_cZ_D{}^d C_{ab}{}^c{}_d-2X^{C}Z_{D}{}^d\nd_{[a}\Rho_{b]d}
- 2X_{D}Z^{Cd}\nd_{[a}\Rho_{b]d}.
$$
The upper case tractor indices may be raised and lowered using the
tractor metric, and the tractor projectors $Z$ and $X$ combine via
the tractor metric according to $Z_D{}^d Z^D{}_e=\delta^d_e$ (the
Kronecker delta), $X^D Z_D{}^d =0$ and $X_DX^D=0$. It follows that
all pairings vanish other than $ZZ$ pairings, and the first
equality in \nn{tCR} is then immediate.

The equality relating $s_{2k}(C)$
to $s_{2k}(R)$ follows easily from a direct calculation; see e.g.\
\cite{Avez}. \qquad $\Box$

\section{The $Q_k$ operators}
In \cite{ddd} the authors constructed, for $k=0,1,\cdots ,n/2+1$,
natural operators $Q^g_k: \ce^k\to
\ce_k$, with $Q^g_{n/2}$ a nonzero constant, $Q^g_{n/2+1}=0$ and
otherwise with properties as follows.
(Here and below, $\cE^k$ denotes the bundle of $k$-forms or,
by way of notational abuse, the sections of this bundle; and
$\cE_k$ is the tensor product of $\cE^k$ with the $(2k-n)$-densities.)
We shall sometimes suppress the dependence of $Q^g_k$ on the metric $g$
and write simply $Q_k$.

\begin{proposition}\label{Qks} Up to a non-zero constant scale,
$Q^g_k$ has the form
$$
(d\d)^{n/2-k}+\LOT
$$
and $Q_01$ is the \IT{Branson}
$Q$-curvature.  Upon restriction to the closed $k$-forms $\cC^k$, $Q_k^{g}$
has the conformal transformation law
\begin{equation}\label{Qoptrans}
Q_k^{\wh{g}}u=Q_k^g u +  \d Q_{k+1}^g d (\om u)
\end{equation}
where $ \wh{g}=e^{2\om}g$ with $ \om$ a smooth function.
\end{proposition}
\noindent It follows that
$\d Q_{k+1}^g d= :L_k:\ce^k\to\ce_k$ is conformally invariant.

A question of some interest is as follows. If $S$ is a natural
Riemannian density with $\int S$ conformally invariant on compact
manifolds, then is it necessary that
$$
S={\rm const} \cdot {\sf Pff}+L+V,
$$
where $L$ is a local conformal invariant and $V$ is an exact
divergence?

This question has implications for the Q-curvature (and other natural
densities in the linear Q-space), since $\int Q$ is conformally
invariant.  From our current viewpoint, $Q$ arises as $Q_0 1$, where
$Q_0$ is a natural differential operator as in Proposition \ref{Qks}.
Rephrasing things to respect this viewpoint, we come to the following
equivalent question. Let $\cC_k:=\cE_k/\cR(\d)$, and let
$\pi_k:\cE_k\to\cC_k$ be the quotient map. Suppose $S$ is a natural
differential operator $\ce^0\to \ce_0$ with the property that
$\pi_0\circe S|_{\cC^0}$ is conformally invariant.  Is it necessary
that
$$
S|_{C^0}=[{\rm const} \cdot {\sf Pff}+L+ \d U]|_{\cC^0},
$$
where $L:\ce^0\to\ce_0$ is a natural differential operator with
$L|_{\cC^0}$ conformally invariant, and $U$ is a natural
differential operator $U:\ce^0\to\ce_1\,$? (Note that the
composition just mentioned acts $\d U:\ce^0\to\ce_0\,$.)
Of course in the last display we view ${\rm
const} \cdot {\sf Pff}$ as a multiplication operator.

More generally, there is an analogue of this question for operators on
$k$-forms with $k=0,1,\cdots ,n/2-1$.  If $S:\ce^k\to \ce_k$ is a
natural differential operator with $\pi_k\circe
S|_{\cC^k}:\cC^k\to\cC_k$ conformally invariant, is it necessary that
\begin{equation}\label{formS}
S|_{\cC^k}=[ {\sf O}+L+ \d U]|_{\cC^k},
\end{equation}
where $L:\ce^k\to\ce_k$ is a natural differential operator with
$L|_{\cC^k}$ conformally invariant, $U$ is a natural differential
operator $\ce^k\to\ce_{k+1}$, and ${\sf O}$ means ${\rm const}
\cdot {\sf Pff}$ for $k=0$, and $0$ for all other $k$?

It turns out that that we can answer this more general question in
the negative if $k\geq 1$, with the $Q_k$ operators as
counterexamples.
Note that by \nn{Qoptrans}, the invariance hypothesis on
$\pi_k\circe Q_k|_{\cC^k}$ is satisfied.
Let us fix $k\in \{1,\cdots ,n/2-1\}$ and
suppose, with a view to contradiction, that as an operator on
$\cC^k$ we have $Q_k =L+ \d U $ with $L$ and $U$ as above. Now we
consider the situation on a conformally flat manifold. Since
$$
L:\cC^k\to \ce_k
$$
is conformally invariant, it follows that
$$
L d:\ce^{k-1}\to \ce_k
$$
is conformally invariant. From the classification of invariant
operators on locally conformally flat manifolds (see e.g.\
\cite{ESlo}) it follows that $Ld=0$. Thus $\d Q_k d=0$ on
conformally flat manifolds. This contradicts
Proposition \ref{Qks}, which asserts that
the leading term of $\d Q_k d$ is a
non-zero constant multiple of $(\d d)^{n/2+1-k}$. In summary we
have the following:

\begin{theorem}\label{nogo}
For $k\in \{1,\cdots n/2-1\}$ the operators $Q_k:\ce^k\to\ce_k$ of
Proposition \ref{Qks} cannot be written in the form $L+\d U$ where
$L:\ce^k\to\ce_k$ is a natural differential operator with
$L|_{\cC^k}$ conformally invariant, and $U$ is a natural
differential operator carrying $\ce^k\to\ce_{k+1}$.
\end{theorem}

Without speculating on any possible implications that the $Q_k$
might have for topology,
this theorem says that from the viewpoint of conformal geometry,
and modulo local conformally invariant operators and divergence
type operators,
the $Q_k$ are $k$-form analogues
of the Pfaffian.  That is, for forms of degree $k\ge 1$, they play
at least one of the roles that the Pfaffian plays for $k=0$.

\section{Conformal forms and a period map}\label{Theta}

Suppose now that $M^{n~{\rm even}}$ is compact without boundary (but
not necessarily oriented or connected) and equipped with a Riemannian
signature conformal structure $[g]$. We construct a family of
invariant bilinear forms.

Recall that a conformal structure $[g]$ is equivalent to a canonical
non-degenerate symmetric bilinear form $\bg$ that takes values in
densities of weight 2. We term this the {\em conformal metric} (see
e.g.\ \cite{ddd} for further discussion) and if $g$ is a metric in the
conformal class then $g=\si^{-2}\bg$ for some non-vanishing weight 1
density $\si$.  For $k$-forms $\xi,\eta$, let $\xi\cdot\eta$ be the
local form inner product determined by $\bg^{-1}$, and
$\langle\xi,\eta\rangle:=\int\xi\cdot\eta$.  Consider the quadratic
form
$$
\cQ_k:\cC^k\times\cC^k\to\bbR,\qquad\cQ(\xi,\eta)=\langle\xi,Q_k\eta\rangle.
$$
Since the operator $Q_k$ is formally self-adjoint at any metric
from $[g]$, the form $\cQ_k$ is symmetric.

Let $\cH_G^k:=\{\xi\in\cC^k\mid \d Q_k\xi=0\}$.  Because the $Q_k$
acts as $(d\d)^{n/2-k}+\LOT$ on closed forms, the system of
equations (at each given conformal scale) satisfied by elements of
$\cH^k_G$ is elliptically coercive, and thus $\cH_G^k$ is finite
dimensional. In more detail, each $\xi\in\cH^k_G$ satisfies a
fixed system of the form
$$
\left(\begin{array}{c}d \\
\d\left\{(d\d)^{n/2-k}+\LOT\right\}\end{array}\right)\xi=0,
$$
and so also satisfies the Laplace-type equation
$$
\begin{array}{rl}
0&=\left((d\d)^{n/2-k}\d\ \mid\  d\right)\left(\begin{array}{c}d \\
\d\left\{(d\d)^{n/2-k}+\LOT\right\}\end{array}\right)\xi \\
&=\left\{(\d d+d\d)^{n/2-k+1}+\LOT\right\}\xi. \end{array}
$$
Let $\tilde\Theta_k$ be the restriction of $\cQ_k$ to a quadratic
form $\cH^k_G\times\cH^k_G\to\bbR$.

Let $H^k_G$ be the image of $\cH^k_G$ in the cohomology $H^k$,
under the cohomology class map $\xi\mapsto[\xi]$. 
Note that is a conformally invariant subspace of $H^k$. 
 We claim that
$\tilde\Theta_k$ gives rise to a quadratic form $\Theta_k$ on
$H^k_G$ via
$$
\Theta_k([\xi],[\eta])=\tilde\Theta_k(\xi,\eta).
$$
To see that this is well-defined, let $\xi,\xi'\in\cH_G^k$ with
$\xi'-\xi=df$.  Since
$$
\tilde\Theta_k(df,\eta)=\langle df,Q_k\eta\rangle=\langle
f,\underbrace{\delta Q_k \eta}_{=0}\rangle =0,
$$
we have $\tilde\Theta_k(\xi',\eta)=\tilde\Theta_k(\xi,\eta)$.

In \cite{ddd}, a conformal manifold $(M,[g])$ is called
$(k-1)${\em-regular} if the map $\cH_G^k\to H^k$ above is
surjective; that is, if $H^k_G=H^k$.  (We expect this condition to
hold for generic conformal manifolds $(M,[g])$.)

\begin{theorem}\label{capform}
If $(M,[g])$ is $(k-1)$-regular, then $\Theta_k:H^k\times
H^k\to\bbR$ is a conformally invariant quadratic form on
cohomology.  $Q_{n/2}$ is the identity, every $(M,[g])$ is
$(n/2-1)$-regular, and $\Theta_{n/2}$ is the integrated metric 
pairing.
Every $(M,[g])$ is $(-1)$-regular, and in the conformally flat
case, up to a universal positive constant multiple, $\Theta_0=
{\rm diag}(\chi(M_1),\cdots,\chi(M_p))$, where
$M=M_1\sqcup\cdots\sqcup M_k$ is the connected component
decomposition of $M$.
\end{theorem}

\noindent{\em Proof}: Everything is clear or established above,
except the final statement.  By \cite{bgp}, in the conformally
flat case, each natural density with conformally invariant
integral takes the form $c\cdot{\rm Pff}+\d\r$, where $c$ is a
constant and $\r$ is a natural one-form-$(2-n)$-density. Applying
this to $Q=Q_01$ and integrating, we find that
$$
c\cdot\chi(M_j)=\int_{M_j}Q.
$$
The round sphere, where $Q=(n-1)!$, serves to normalise $c$:
$$
2c=(n-1)!\varpi_n\,,
$$
where $\varpi_n$ is the volume of the round sphere. Each connected
component contributes one dimension of 0-cohomology (generated by
the class of its characteristic function); this establishes that
$$
\Theta_0=\dfrac{(n-1)!\varpi_n}2{\rm
diag}(\chi(M_1),\cdots,\chi(M_p)).\qquad\Box
$$
\vspace{0.2cm}

The theorem indicates that in some sense, the conformally invariant
forms $\Theta_k$ interpolate between the integrated metric pairing
of middle-forms and the Euler characteristic.

In another direction, we get a map from the set of conformal structures
that generalises the celebrated period map, but involves forms
not of the middle degree.
Recall that above we used $\cE_k$ to denote the tensor product of
$\cE^k$ with the $(2k-n)$-densities. Given a conformal structure we
may use the conformal metric to canonically identify this space with
the tensor product of $\Lambda^k TM $ with $(-n)$-densities. Then we
again use $\d$ to denote the adjoint of $d$ afforded by the global
pairing of $\ce^k$ with $\ce_k$ and we write $H_k$ for the
corresponding cohomology space. On a fixed conformal structure this
change makes no difference at all. On the other hand this
contravariant point of view has the advantage that the operators $\d$
and the cohomology spaces $H_k$ are now just objects belonging to the
smooth structure, and in particular are not affected by moving the
conformal structure.

For $0\le k\le n/2$, let 
us equip $H^k\oplus H_k$ with the obvious symplectic structure. 
We obtain the following generalisation of the
period map, modelled on the construction of Eastwood and Singer
in \cite{EastSin2} for the case $n=4$, $k=1$.
\begin{proposition}\label{periodprop}
For each conformal structure $[g]$ on $M$, there is a well-defined map 
$$
I^{[g]}:\cH^k_G \to H^k\oplus H_k
$$ given by $\phi\mapsto ([\phi],[Q_k \phi])$, the
range of which  
is a Lagrangian subspace. Thus we obtain a map
$$
\Phi:\{\mbox{conformal structures on M}\}\to G_{M}
$$
where $G_M$ is the Grassmannian of Lagrangian subspaces of the 
symplectic vector space $H^k\oplus H_k$, given by $\Phi([g])=\cR(I^{[g]})$.
\end{proposition}
\noindent{\bf Proof:} 
By the remarks preceding the proposition, we may regard $H^k\oplus H_k$
as a fixed target space, independent of the conformal structure.
{}From the transformation law \nn{Qoptrans} it
follows that the map $\phi\mapsto ([\phi],[Q_k \phi])$ depends only on
conformal structure.  (Or see section 2.0 of \cite{ddd} where
$I^{[g]}$ is also discussed).  That $\cR(I^{[g]})$ is an isotropic
subspace follows immediately from the result \cite{ddd} that, in any
scale $g$, $Q^g_k$ is formally self-adjoint. Then $\dim(\cR(I^{[g]})$
is the $k^{\rm th}$ Betti number from Corollary 2.12 of \cite{ddd} and
so $\cR(I^{[g]})$ is Lagrangian (i.e.\ it is isotropic and has half the
dimension of $H^k\oplus H_k $).  $\Box$

\section{Natural scalars in the linear Q-space}

For a metric $g$ from the conformal class and $k\in \{0,1,\cdots
,n/2+1 \}$, let us fix a pair $(\xi,\eta)$ in $\cC^k\times \cC^k$
and make the definitions
\begin{equation}\label{newQdef}
\tilde Q^g_{\xi,\eta}:=\xi\cdot Q_k^g\eta,\qquad
Q^g_{\xi,\eta}:=\frac12(\tilde Q^g_{\xi,\eta}+\tilde
Q^g_{\eta,\xi}), \qquad Q_\xi^g:=Q^g_{\xi,\xi}\,.
\end{equation}
(Recall that $\xi\cdot\eta$ is the pointwise form inner product
determined by the conformal metric $\bg$.) 
This makes each quantity in the display a
density. 
At the moment, we make no naturality assumption on $\xi$
or $\eta$. Related to these constructions are the differential
operators $\tilde L_{\xi,\eta}$, $L_{\xi,\eta}$, and $L_\xi$
carrying $\ce^0$ to $\ce_0$, defined by 
\begin{equation}\label{Ldef}
\tilde L_{\xi,\eta} f: = \xi\cdot L_k (f \eta),\qquad
L_{\xi,\eta}=\frac12(\tilde L_{\xi,\eta}+\tilde
L_{\eta,\xi}),\qquad L_\xi=L_{\xi,\xi}\,.
\end{equation}
Note that by construction each operator in \nn{Ldef} is conformally
invariant.
Some properties of these  are summarised in the following proposition.

\begin{proposition}\label{newQprop}
\IT{i} At a fixed conformal scale $g$, $Q_{\xi,\eta}^g$ is a
density. Under a conformal change of metric to $\wh g= e^{2\om}
g$ this has the conformal transformation
\begin{equation}\label{Qxetran}
Q_{\xi,\eta}^{\wh g}=Q_{\xi,\eta}^g+L_{\xi,\eta} \om .
\end{equation}
\IT{ii} The operator $L_{\xi,\eta}$ has the form $\d M^g_{\xi,\eta}
d$ and is formally self-adjoint.\\
\IT{iii} On compact manifolds $\int Q_{\xi,\eta}$ is conformally invariant.\\
\IT{iv}
In the case $k=0$ we have that $Q^g_{1,1}$
is the Q-curvature.\\
\end{proposition}

\noindent{\bf Proof:} Part \IT{i}. Since $\tilde Q^g_k$ takes
values in $\ce_k$ and the form inner product carries a density
weight of $-2k$, the scalar densities $\tilde Q_{\xi,\eta}$ have
weight $-n$.  Thus this also holds for $Q_{\xi,\eta}$.
The transformation law is an immediate
consequence of the transformation law \nn{Qoptrans}  for
the $ Q_{k}^g$ operator.
(In fact, the stronger statement obtained by attaching a tilde to each
quantity in \nn{Qxetran} also holds.)

Part \IT{ii}.
Since $\eta$ is closed,
$$
L_k (f_1 \eta)= \d Q_{k+1}^g \big( \e(d f_1) \eta\big)
$$
for any function $f_1$. Thus $\tilde{L}_{\xi,\eta}$
has $d$ as a right composition factor. On the other hand for another
function $f_2$, from the fact that $\xi$ is closed we have
$$
\langle
f_2 \xi , \d Q_{k+1}^g \big( \e(d f_1) \eta\big) \rangle= \langle \e(d
f_2) \xi , Q_{k+1}^g \big( \e(d f_1) \eta\big) \rangle,
$$
where recall $\langle\phi,\psi\rangle=\int\phi\cdot\psi$.
Thus $\d$ is a
left composition factor of $\tilde{L}_{\xi,\eta}$ and overall we have
that $\tilde{L}_{\xi,\eta} = \d \tilde{M}^g_{\xi,\eta} d$ for some
differential operator $\tilde{M}^g_{\xi,\eta} $. Continuing our
integration by parts we see that
\begin{equation*}\begin{split}
\langle d (f_2 \xi) , Q_{k+1}^g d (f_1 \eta) \rangle
& = \langle Q_{k+1}^g d (f_2 \xi) , d (f_1 \eta ) \rangle
= \langle d (f_1 \eta), Q_{k+1}^g  d(f_2 \xi) \rangle\\
&= \langle f_1 ,  \eta \cdot \d Q_{k+1}^g d ( f_2 \xi) \rangle
\end{split}\end{equation*}
and so $\tilde{L}_{\eta,\xi}$ is the formal adjoint of
$\tilde{L}_{\xi,\eta}$.
Symmetrising over $\xi$ and $\eta$ the
claimed results follow.

Part  \IT{iii} is immediate from parts \IT{i} and \IT{ii}.
Finally, by definition $Q^g_{1,1}=Q_0^g 1 $ and so \IT{iv} is
immediate from the result $Q=Q_0^g 1  $ in Proposition \ref{Qks}.
\quad $\Box$
\vspace{0.2cm}

An interesting possibility is to take $\xi,\eta$ in \nn{newQdef}
to be invariant forms arising from a connection $D$, for example
$\xi=\eta = s_{\ell}(F^D)$. Then (for example) $\int Q_{\xi}^g$
gives some conformal coupling between the geometric structure and
the connection $D$.

In the case that $D$ is the conformal tractor connection this
yields {\em natural} invariants $Q_{\xi}$.
As a result, we have
the following.

\begin{theorem}\label{naturalQs}
  For each pair of conformal $k$-forms $\tau, \kappa$, where $0\leq
  4k\leq n/2-1$, the scalar $Q^g_{\tau,\kappa}$ field is a natural
  invariant in the linear Q-space $\cIQ$. On compact manifolds $\int
  Q^g_{\tau,\kappa}$ is an invariant of the conformal structure.
\end{theorem}

\subsection{A coupled generalisation of $Q_{n/2-1}$}

The operator
$$
Q_{n/2-1}:\ce^{n/2-1}\to \ce_{n/2-1}
$$
has the explicit
formula $Q^g_{n/2-1}=d \d -4 \Rho\sharp + 2\J$. Here we view the
Schouten tensor $\Rho$ as an endomorphism of the tangent bundle and
$\sharp$ indicates the usual ($\otimes$-derivation)
action of such an endomorphism on tensors
(in this case on forms). $\J$ is the trace of $\Rho$.

Acting
on the closed forms $\cC^{n/2-1}$, we have the conformal transformation
law
$$
Q^{\widehat{g}}_{n/2-1} \kappa = Q^g_{n/2-1} \kappa +2 \d d (\om \kappa) .
$$ This transformation law is preserved if we couple to a connection.
Suppose $D$ is a connection on some vector bundle $\cV$. Then write $
Q^{D,g}_{n/2-1} $ for the operator on the $\cV$-valued $(n/2-1)$-forms
$\ce^{n/2-1}(\cV)$ given by the formula
$$
d^D \d^D -4 \Rho\sharp + 2\J .
$$
Write $\cC^k(\cV)$ for the space of $\cV$-valued $k$-forms $\kappa$
satisfying the identity $d^D\kappa=0$.
By direct calculation one readily verifies the following result.

\begin{theorem}\label{couplethm} For any
vector bundle with connection $(\cV,D)$$$
Q^{D,\widehat{g}}_{n/2-1} \kappa =
Q^{D,g}_{n/2-1} \kappa + 2 \d^D d^D (\om \kappa) ,
$$
for $\kappa\in \cC^{n/2-1}(\cV)$ and
$\widehat{g}=e^{2\om}g$.
\end{theorem}

Since the curvature $F$ of the connection $D$ and its exterior
powers are annihilated by $d^D$, it is straightforward to use these to
construct coupled quantities in the linear Q-space. Once again natural
scalar fields are obtained when we specialise to the case in which $D$ is
the tractor connection $\nabla$ on the standard tractor bundle $\bbT$
and its tensor powers. Writing $\Om$ for the curvature of the
tractor connection on $\bbT$, in dimension 6 we have
\begin{equation}\label{6Dquant}
\Omega^A{}_B \cdot Q^\nd_{2}\Omega^B{}_A
\end{equation}
for example.
Here we have displayed abstract tractor indices but omitted
form indices.
Direct computation yields the following description of this invariant in
terms of classical curvatures.  Let $A_{abc}:=2\nd_{[c}\V_{b]a}$ be the
Cotton tensor, let
$$
W_{abcde}:=\nd_eC_{abcd}+2g_{e[a}A_{b]cd}+2g_{e[c}A_{d]ab},
$$
and let
$$
U_{abcd}:=\nd_aA_{bcd}-\V_a{}^eC_{ebcd}.
$$
Then \cite{fgast,crgsrni99,CapGoamb}
$$
I:=|W|^2-16(C,U)+16|A|^2
$$ is conformally invariant in dimension 6.  Gover and Peterson
\cite{GoPetLap} note that $G:=\Delta|C|^2=\delta d|C|^2$ admits the
linear and formally self-adjoint conformal change law
$$
\wh{G}=G+4\d(|C|^2d\om).
$$
Fefferman and Hirachi \cite{FHmrl} note that
$$
H:=-C_{abcd}C^{abce}\V^d{}_e+|A|^2+\frac14|C|^2\J
$$
has a linear conformal change law, which can be shown by direct computation to
be formally self adjoint; in fact,
$$
\wh{H}=H+\left(\frac14|C|^2\Delta+4\nd_c\V_{ab}C^{acbd}\nd_d+C^a{}_{cde}
C^{bcde}\nd_a\nd_b\right)\om.
$$
The invariant \nn{6Dquant} is
$$
\Omega^A{}_B \cdot Q^\nd_{2}\Omega^B{}_A=\frac14I
+\frac18G-2H-
\frac14C_{abcd}C^{ab}{}_{ef}C^{cdef}-C_{abcd}C^a{}_e{}^c{}_fC^{bedf}\,.
$$

In dimension 10, examples include (skewing over the unexpressed tensor
indices in each $\Om\otimes \Om$)
$$
\Omega^A{}_B \Omega^B{}_C \cdot Q^\nd_{4}\Omega^C{}_E \Omega^E{}_A  \quad
\mbox{ and }
\quad \Omega^A{}_B \Omega^C{}_E \cdot Q^\nd_{4}\Omega^E{}_C \Omega^B{}_A .
$$

One might expect that there is a result generalising Theorem \ref{couplethm}
to the $Q_{k}$ for $k\leq n/2-2$. This is not the case.
\begin{theorem}\label{nocouplethm}
Suppose that
$$
Q^{D,g}_{n/2-k}:\cC^k(\cV)\to\ce_k(\cV), \qquad k\in\{
0,1,\cdots n/2-2\}
$$
is given by a universal polynomial formula, with natural coefficients, 
in the connection $D$
coupled with the Levi-Civita connection.
Suppose further that for any vector
bundle/connection pair $(\cV,D)$ we have
\begin{equation}\label{multop}
Q^{D,\widehat{g}}_{n/2-k} = Q^{D,g}_{n/2-k} + S^{g,D} d^D \om
\end{equation}
whenever $\widehat{g}=e^{2\om}g$.
{\rm(}In the last term of {\rm(}\ref{multop}{\rm)} we
view $\om$ as a
multiplication operator.{\rm)}
Then the operator $S^{g,D} d^D$ vanishes on conformally flat structures.
\end{theorem}

\noindent{\bf Proof:}
It is straightforward to show that without loss of generality we
may assume that $S^{g,D}$ is also given by a universal polynomial
formula, with natural coefficients, 
in the coupled connection.
{}From formula \nn{multop} it follows that $S^{g,D}d^D$ is conformally
 invariant acting on $\w\cC^k(\cV)$ for each function $\w$.  Thus by
 linearity we obtain an operator $S^{g,D} d^D:\ce^k(\cV)\to
 \ce_k(\cV)$ which is conformally invariant for any vector
 bundle/connection pair $(\cV,D)$. But Proposition 1.1 of
 \cite{Gosrni03} shows that any such operator vanishes on conformally
 flat structures.  \quad $\Box$

\subsection{A more complete picture}\label{complete}
Although it is a slight digression we should point out that many of
our constructions generalise without difficulty in an obvious way. For
example Proposition \ref{newQprop} shows that $Q_k$ generates scalar
fields in the linear Q-space via \nn{newQdef}. However if we take
$(\xi,\eta)\in \cC^{k-\ell}\times \cC^{k-\ell}$, for $\ell\leq k$,
then the formula \nn{newQdef} may be used to instead give an operator
$Q_{\xi,\eta}^{k,g}:\ce^{\ell}\to \ce_{\ell}$.
In abstract index notation, we take (up to nonzero constant multiples)
$\xi^{a_{k-\ell+1}\cdots a_k}
Q_{a_1\cdots a_k}{}^{b_1\cdots b_k}\eta_{b_{k-\ell+1}\cdots b_k}$
to obtain $(\tilde Q_{\xi,\eta})_{a_1\cdots a_\ell}{}^
{b_1\cdots b_\ell}$, and set $Q_{\xi,\eta}=\frac12
(\tilde Q_{\xi,\eta}+\tilde Q_{\eta,\xi})$.
These maps have properties as
follows.

\begin{proposition}\label{newQpropforms}
\IT{i} $Q^{k,g}_{\xi,\eta}$ is a differential operator which,
upon restriction to the closed forms
$\cC^{\ell}$, has the conformal transformation law
$$
Q_{\xi,\eta}^{k,\wh g}=Q_{\xi,\eta}^{k,g}+  \d Q_{\xi,\eta}^{k+1,g}d \om \quad \mbox{where} \quad
\wh g= e^{2\om} g,
$$
 and on the right hand side $\om$ is viewed as a multiplication operator.\\
\IT{ii} The operator  $\d Q_{\xi,\eta}^{k+1,g}d$ is formally self-adjoint. \\
\IT{iii} $ Q_{\xi,\eta}^{k,g}$ determines a conformally invariant operator
$$
 Q_{\xi,\eta}^{k}:\cC^\ell\to \cC_\ell:=\cE_\ell/\d\cE_{\ell+1}\,.
$$ 
\IT{iv} On compact manifolds $ Q_{\xi,\eta}^{k,g}$ induces a
conformally invariant pairing between ${\cN}( \d Q_{\xi,\eta}^{k+1,g}d)$ and $\cC^\ell$ given by
$$
(\mu,\rho)\mapsto \int \mu\cdot Q_{\xi,\eta}^{k,g} \rho .
$$
In particular we get a pairing on $\cC^\ell\times \cC^\ell$.\\
\IT{v} $\d Q_{\xi,\eta}^{k,\wh g}$ is conformally invariant on the
null space of $\d Q_{\xi,\eta}^{k+1,g}d:\ce^\ell\to \ce_\ell$. Thus in
particular it is conformally invariant on $\cC^\ell$.
\end{proposition}

\noindent{\bf Proof:} Parts \IT{i} and \IT{ii} are proved by
obvious adaptations of the earlier arguments. Parts \IT{iii} and
\IT{iv} then follow immediately. It is shown in \cite{ddd} that on
$\ce^k$ we have $\d Q_k^{\wh g}= \d Q_k^{\wh g} +c\i (d\w) \d
Q_{k+1}^g d$, for some constant $c$.
Part \IT{v} follows easily. \quad $\Box$\\

The Proposition shows that the operators $Q_{\xi,\eta}^{k,g}$
generalise in a natural way the operators $Q_k^g$ of \cite{ddd}.
Many of the other results for the $Q_k^g$ carry over to the
$Q_{\xi,\eta}^{k,g}$. In fact the story is still more general
since we can compose the $Q_k$ on the right by exterior
multiplication with any closed form  and on the left by interior
multiplication with any closed form.  In these matters we are making
no attempt be complete in the current note.

By using conformal forms for the $\xi,\eta$ in these constructions the
operators in the proposition become natural operators. Of course since
the conformal forms vanish on conformally flat structures the
resulting operators all vanish (apart from the $Q_k^g$), and this includes
the natural
scalar fields in the
linear Q-space that we constructed earlier as a special case.

\section{An example}\label{cpH}

Consider the case of a 10-dimensional manifold $M=\cph\times N$,
where $N$ is a conformally flat 6-dimensional manifold of constant
scalar curvature $\nu$. The factor $\cph$ is supplied with the
Fubini-Study metric. Under these circumstances the Pontrjagin
4-form $\k$ of $M$ is just (i.e.\ pulls back to under inclusion
of a submanifold $\cph\cong\cph\times\{y\}$) 
that of $\cph$, and this is an
eigenform of $Q_4$ with eigenvalue a nonconstant affine function
of $\nu$.  Thus, except for one value of $\nu$, the quantity
$Q_\k$ is nonzero.

First we claim that on general oriented Riemannian 4-manifolds,
the Pontrjagin 4-form of the Levi-Civita connection is
\begin{equation}\label{pontclaim}
p_1^{{\rm LC}}=\dfrac1{96\pi^2}(|C_+|^2-|C_-|^2)E,
\end{equation}
where $E$ is the volume form.  
In fact, from 
\nn{pontdef} and Proposition \ref{confformeq}, we have
$$
p_1^{{\rm LC}}=-\frac1{8\pi^2}s_2(C).
$$
An elementary calculation gives
$$
\begin{array}{rl}
4!g(s_2(C),E)&=C^a{}_{bcd}C^b{}_{aef}E^{cdef}=2C^a{}_{bcd}(\star C)^b{}_a{}^{cd} \\
&=2(C_++C_-)^a{}_{bcd}(C_+-C_-)^b{}_a{}^{cd}=-2(|C_+|^2-|C_-|^2),
\end{array}
$$
so that
$$
s_2(C)=-\frac1{12}(|C_+|^2-|C_-|^2)E.
$$
Equation \nn{pontclaim} follows immediately.

Now consider
the special case of $\bbC\bbP^2$ with the Fubini-Study metric
$g$. $\bbC\bbP^2$ is orientable, its Euler characteristic is
$\chi=3$, and its signature is $\si=1$.  The usual metric
normalisation has
scalar curvature 24 (and thus Schouten scalar 4) and volume
$\pi^2/2$. The Pfaffian in dimension 4 is
$$
{\rm Pff}=\dfrac1{32\pi^2}(|C|^2-8|S|^2+6J^2),
$$
where $S$ is the trace-free Schouten tensor.
Since $g$ is locally symmetric, $\nd C$ vanishes and $|C|^2$ is constant.
Since $g$ is Einstein, $S=0$; this and the above data give
$|C|^2=96$. The Hirzebruch polynomial (signature
integrand) in dimension 4 is
$$
L=\dfrac1{48\pi^2}(|C_+|^2-|C_-|^2),
$$
so $|C_+|^2=96$ and $|C_-|^2=0$, reflecting the fact that
$\bbC\bbP^2$ is half conformally flat, and indicating which half
is flat. In particular,
$$
p_1^{{\rm LC}}(\bbC\bbP^2,g)=E/\pi^2,
$$
where $E$ is the volume form of $g$, and
$\int_{{\bbC\bbP}^2}p_1^{{\rm LC}}=\frac12$.

Now consider the 10-dimensional manifold $\cph\times N$.  Since the
Riemann curvature tensor is additive over products it follows easliy
that  the Pontrjagin forms of a direct
product depend only on the Weyl tensors of the factors.  
Now
$$
Q_4\k=(d\d-4\Rho\sharp+2\J)\k=2(\J-2\Rho\sharp)\k,
$$
since $\k$ is (parallel and hence) harmonic. The scalar curvature of $M$ is $24+\nu$,
so the Schouten scalar is $\J=(24+\nu)/18$.  
The pullback to $\cph$ of the Schouten tensor $\Rho$ of $M$ is 
$(84-\nu)g/144$, so $\Rho\sharp\kappa=(84-\nu)\kappa/36$.
Thus $\k$ solves the eigenvalue problem
$$
Q_4\k=2(\nu-30)\k/9
$$
(see Section \ref{moments}).
Since by construction $\nd$ annihilates $Q_\k=\k\cdot Q_4\k$,
the function corresponding to this density under the metric $g$
is constant, and the above shows that this constant  
is nonzero as long as $\nu\ne
30$.  In particular, $\int Q_\kappa$ is nonzero if $\nu\ne 30$.

\section{Invariant nonlinear functionals and prescription problems} \label{pres}

For the moment, suppose our conformal manifold is compact, but not
necessarily of Riemannian signature.
For each quantity $Q$ in the linear Q-space, we define a two point
functional $\cK(\wh g,g)=\cK_{Q}(\wh g,g)$ on the conformal class
$[g]$ by
$$
\cK(\wh g,g)= \tfrac12\int_M \om(\wh g,g) (Q^g+Q^{\wh g}),
$$
where $\om(g,\wh g)$ is the unique function on $M$ satisfying $\wh
g=e^{2 \om (g,\wh g)}g$.  As a $\ci(M)$-valued function on the
conformal class $[g]$, $\w$ is a {\em cocyle}:
$$
\om(g_2,g_1)=-\om(g_1 , g_2),\qquad
\om(g_3,g_2)+\om(g_2,g_1)=\om(g_3,g_1)
$$
for any metrics $g_1,g_2,g_3\in[g]$. The real valued function
$\cK$ is also a cocycle in this sense: it is alternating because
$\w$ is, and with $\w_{ij}:=\w(g_i,g_j)$ and $Q_j:=Q^{g_j}$, we have 
$$
\om_{21} (Q_1+Q_2) + \om_{32} (Q_2+Q_3)
=2\w_{31}Q_1+\w_{21}L\w_{21}+\w_{32}L\w_{32}+2\w_{32}L\w_{21}\,,
$$
using the conformal invariance of $L$ and the conformal
transformation law for $Q$.
This differs from
$$
\om_{31} (Q_1+Q_3) = 2 \om_{31}Q_1+ (\om_{21}+\om_{32}) L
(\om_{21}+\om_{32})
$$
by
$$
\om_{32} L \om_{21}-\om_{21} L \om_{32}\,,
$$
which is a divergence since $L$ is formally self-adjoint. This
proves that
$$
\cK(g_3, g_2)+ \cK(g_2, g_1)= \cK(g_3,g_1).
$$

Fixing a metric $g_0$, we wish to look for metrics that are
critical for the functional $\cK(g,g_0)$ with respect to conformal
variations of $g$.  Because uniform scaling of $g$ results in the
addition of constant multiples of the conformal invariant $c:=\int
Q$ to the functional, we can either restrict ourselves to
equal-volume perturbations (thus freezing out uniform scaling), or
add a volume penalty to the functional, as in:
$$
\cM(g,g_0)=-\dfrac{c}{n}\log\dfrac{{\rm vol}(g)}{{\rm vol}(g_0)}
+\cK(g,g_0)=:\cV(g,g_0)+\cK(g,g_0).
$$
Note that the functional $\cV$ depends on our Q-quantity through
the constant $c$.  The functional $\cM(g,g_0)$ is invariant under
uniform scaling: addition of a constant $b$ to $\w(g,g_0)$
increases $\cK$ by $bc$, and decreases $\cV$ by the same amount.
In addition, the penalised functional $\cM$ is a cocycle on $[g]$.
Thus the problem of finding critical metrics or of extremising
$\cM(g,g_0)$ is independent of the choice of base metric $g_0$,
since choosing base metric $g_1$ instead just adds the constant
$\cM(g_0,g_1)$ to the functional.
\vspace{.4cm}

\noindent{\bf Remark}: In what follows, we shall integrate both
densities and functions.  The integral of a density, for example
$\int Q$, makes sense given only a conformal class.  To integrate
a function $f$, we need to choose a metric $g$ within the conformal
class and use its (pseudo-)Riemannian measure; we shall
denote this process by $\int f\,dv_g$.  Though 
$Q$ is a density, the notion of a
metric $g$ of {\em constant} $Q$ is well-defined.  The choice $g$
of scale induces an identification of all the density bundles
$\cE[w]$; in particular there is a canonical function ($0$-density)
corresponding to $Q$ via the metric $g$, to which we also give the
name $Q$.  
\vspace{.4cm}

Now restrict to the Riemannian signature setting, keeping the
underlying manifold $M$ compact. Taking an equal-volume conformal
curve of metrics $e^{2\b_t}g$, with $\beta_0=0$, and differentiating
to find the variation, we have
$$
\begin{array}{l}
\dfrac{d}{dt}\cM(e^{2\b_t}g,g_0)\bigg|_{t=0}
=\dfrac{d}{dt}\cK(e^{2\b_t}g,g_0)\bigg|_{t=0} 
\\
=\dfrac{d}{dt}\left(\cK(g,g_0)+\cK(e^{2\b_t}g,g)\right)\bigg|_{t=0}
=\displaystyle\int\b Q^gdv_g,
\end{array}
$$ 
where
$\b:=(d/dt)\b_t|_{t=0}\,$.  As $\{\b_t\}$ runs through all
equal-volume conformal perturbations, $\b$ runs through all smooth
functions with $\int\b dv_g=0$; that is, all $\b$ which are $g$-orthogonal 
to the
constants.  The critical metric condition on $g$ is just $Q^g={\rm
const}$.  In fact, the value of the constant at a critical metric $g$
is determined by the conformal invariance of $c=\int Q$ to be
$Q^g=c/{\rm vol}(g)$. Note that because $\cM$ is invariant under
uniform scaling, it was sufficient (in determining the critical
points) to take equal-volume perturbations.
\vspace{.4cm}

\noindent{\bf Remark}: The {\em function} $Q$ 
satisfies the exponential prescription
equation 
\begin{equation}\label{expopresc}
P\om+Q=\wh{Q}e^{n\om},
\end{equation}
where $P$ is the operator on functions obtained from the 
scale-induces density/function correspondence described in the last
remark.
The behavior of this 
prescription equation as a PDE problem varies with the order of
$L$. At the high-order end is the original Q-curvature, whose $L$ has
the form $\Delta^{n/2}+\mbox{LOT}$, and in particular is
positively elliptic. (A special case is Gauss curvature
prescription in dimension 2.) At the other extreme are
Q-quantities with $L=0$; i.e. local conformal invariants.  For
such a quantity, the prescription problem is just algebraic:
Given $g$, we can find a conformally related
$\wh{g}$ with constant $\wh{Q}$ if and only if the sign ($+$, $-$,
or $0$) of $Q$ is constant; and then such a $\wh{g}$ is unique up to
uniform dilations. 
\vspace{.4cm}

Suppose that we have an $\cM$-critical metric $g_0$, i.e.\ a
metric with constant $Q_0$, and consider the corresponding 
one-metric functional
$$
\cH^{g_0}(g):=\cM(g,g_0)=-\dfrac{c}{n}\log\dfrac{{\rm vol}(g)}{{\rm
vol}(g_0)} +\cK(g,g_0),
$$ 
where $c$ is the conformal invariant $\int Q$.

\begin{theorem}\label{invtfcnl}
If $Q_0$ is constant, the functional $\cH^{g_0}(g)$
is invariant under conformal diffeomorphisms $h$ in the identity
component of the conformal group of $(M,[g_0])$, in the sense that
$\cH^{g_0}(h\cdot g)=\cH^{g_0}(g)$, where on covariant tensors
$h\cdot=(h^{-1})^*$.
\end{theorem}
\noindent{\bf Proof}: We need to know that
$$
\cM(h\cdot g,g_0)\stackrel{?}{=}\cM(g,g_0).
$$
We know that
\begin{equation}\label{provis}
\cM(h\cdot g,h\cdot g_0)=\cM(g,g_0),
\end{equation}
since this
is true of the $\cV$ functional (diffeomorphism does not change
the volume), while for the $\cK$ functional,
$$
2\cK(g,g_0)=\int\w(Q_0+Q)=\int h\cdot\left\{\w(Q_0+Q)\right\}
=\int(h\cdot\w)(Q^{h\cdot g_0}+Q^{h\cdot g}),
$$
by naturality of $Q$.  Since $h\cdot\w(g,g_0)=\w(h\cdot g, g_0)$,
this is
$$
2\cK(g,g_0)=\int\w(h\cdot g,h\cdot g_0)(Q^{h\cdot g_0}+Q^{h\cdot
g})=2\cK(h\cdot g,h\cdot g_0).
$$
Using the cocycle condition and \nn{provis}, we have
$$
\begin{array}{rl}
\cH(h\cdot g)&=\cM(h\cdot g,h\cdot g_0)+\cM(h\cdot g_0,g_0)
=\cM(g,g_0)+\cM(h\cdot g_0,g_0) \\
&=\cH(g)+\cM(h\cdot g_0,g_0).
\end{array}
$$
Thus what we need to show is that the very last term $\cM(h\cdot
g_0,g_0)$ vanishes.  Since it vanishes with $\cV$ in place of $\cM$,
what we need to show is the vanishing of
\begin{equation}\label{desideratum}
\begin{array}{rl}
2\cK(h\cdot g_0,g_0)&=\dint\w(h\cdot g_0,g_0)(Q_0+h\cdot Q_0) \\
&=Q_0\dint\w(h\cdot g_0,g_0)(dv_{g_0}+dv_{h\cdot g_0}),
\end{array}
\end{equation}
the last simplification depending on the fact that as a function, $Q_0$ is
constant (and so equals its pushout $h\cdot Q_0$ under $h$).

The corresponding infinitesmimal statement comes from taking
$h=h_t$ in a one-parameter group of conformal transformations, and
computing $(d/dt)|_{t=0}$ of \nn{desideratum}.  This yields
$2Q_0\int\w\,dv_{g_0}\,$, where $\w$ is the infinitesimal conformal factor of
the conformal vector field $T$ generating the $h_t$:
$$
L_Tg_0=2\w g_0\,,
$$
where $L_T$ is the Lie derivative.  
However $n\int\w\,dv_{g_0}$ is the variation of the volume in this
conformal diffeomorphism direction, and so it vanishes.
This establishes invariance of $\cH$ under transformations from
the identity component of the conformal group $G$ of $(M,[g])$ (in
which the one-parameter groups generate a dense set).
$\qquad\Box$
\vspace{0.2cm}

{\bf Remark}: 
For a conformal transformation $h$ in another connected component
$G_1$
of the conformal group of $(M,[g_0])$, we get the same invariance
statement provided $G_1$ contains an isometry of $g_0$.
Indeed, $\cH$ is clearly invariant under
isometries: the integrand of \nn{desideratum} vanishes identically
for the isometry invariance problem.
\vspace{0.2cm}

Given a critical metric $g_0$, an interesting way to rewrite the
functional $\cH(g)=\cM(g,g_0)$ is as follows.  Using the constancy
of $Q_0$ and the conformal transformation law $e^{n\w}Q=Q_0+L_0\w$ (where
$\w:=\w(g,g_0)$), we have
$$
\cH(g)=-\dfrac{c}{n}\log\dfrac{\int e^{n\w}dv_{g_0}}{{\rm
vol}(g_0)}+\left(Q_0\int\w\,dv_{g_0}+\frac12\int\w(L_0\w)dv_{g_0}\right).
$$
Since $Q_0=c/{\rm
vol}(g_0)$, the first term on the right side can be combined with
the left side to give
\begin{equation}\label{towardMT}
\cH(g)=-\dfrac{Q_0{\rm vol}(g_0)}{n}\log\dfrac{\int
e^{n(\w-\overline{\w})}dv_{g_0}}{{\rm vol}(g_0)}+\frac12\int\w(L_0\w)
dv_{g_0}\,,
\end{equation}
where $\overline{\w}$ is the $g_0$-average value of $\w$.

In the case of the original Q-curvature, for which $L$ is the
critical GJMS operator, and of the standard conformal class on the
sphere $S^n$, the quantity on the right hand side of \nn{towardMT} is the
one asserted to be nonnegative by the celebrated
Beckner-Moser-Trudinger inequality \cite{bec}. In this case
$Q=(n-1)!$ (at the round metrics, which constitute one orbit
within the conformal class under the conformal diffeomorphism
group).  What we have done above is to generalise this form and
its conformal transformation invariance to arbitrary Q-quantities:

\begin{corollary}\label{mtcor}
If $Q_0$ is constant, the functional \nnn{towardMT} 
is invariant under conformal diffeomorphisms $h$ in the identity component
of the conformal group of $(M,[g_0])$.
\end{corollary}

\subsection{Generalised prescription problems}
\label{moments}

The prescription problem for the (classical) Q-curvature is
generalised by the prescription/eigenvalue problem
of finding a triple $(\wh{g},\k,\l)\in[g]\times\cH^k_G\times\bbR$ with
the property that $ Q^{\wh{g}}_k\k=\l\k$.  From the viewpoint of an
arbitrary scale $g\in[g]$, this is the differential equation
$L(\w\k)+Q_k\k=e^{(n-2k)\w}\l\k$, where $\w=\w(\wh{g},g)$.  
(We could also state a version of
this problem where we demand only that $\k\in\cC^k$. So one may ask,
for example, whether the problem may be solved with $\k$ a Pontrjagin
form.)  Note that for any solution $(\wh{g},\k,\l)$ with $\l\ne 0$,
the conformal harmonic $\k$ is actually a harmonic at the scale $\wh{g}$, 
since $0=\d Q_k\k=\l\d\k$.  Note also that the example worked out
in Section \ref{cpH} is an example of this eigenvalue phenomenon. If
we nominate and fix a conformal harmonic (or closed form) $\k$, then
the overall problem generalises 
constant Q-curvature prescription: here the general harmonic is a linear
combination of indicator functions for the connected components
of our compact manifold.
In any case where we have a
solution $(\wh{g},\k,\l)$, the conformally invariant quantity
$\langle\k,\wh{Q}\k\rangle$ is the eigenvalue $\lambda$ times $\int
(\k,\k)$, where the latter is computed in the scale $\wh{g}$. In the
usual (i.e.\ $k=0$) problem $\int (\k,\k)$ is a nonnegative
linear combination
of the component $\wh{g}$-volumes, 
in which the coefficients are the squares of those on the indicator functions
in the original expression for $\k$.
Put another
way, we see that the information in $\langle\k,\wh{Q}\k\rangle$
for general form order
generalises that provided in the 0-form case by the conformal invariants $\int_{M_i}Q$
for $M_i$ the connected components of $M$.
The eigenvalue problem described above is a corresponding generalisation of
the problem of precribing constant $\wh{Q}$ on some chosen union of connected
components.

Knowledge of the scale $\wh{g}$ in which the eigenvalue $\l$ is
attained yields further information  on $\l$ and $\k$ which depends explicitly
on $\wh{g}$. To describe this information, note that in the equation
$\langle\k,\wh{Q}\k\rangle=\l\langle\k,\k\rangle$, the left side uses
the conformally invariant quadratic form $\tilde\Theta_k$ (of Sec.\
\ref{Theta}) on $\cH^k_G$, while the right side uses the
(scale-dependent) $\wh{g}$-metric form.  Let $B_k:=\dim\cH_G^k$.
Composing $\tilde\Theta_k$ with the inverse of the $\wh{g}$-form gives
an endomorphism $T^{\wh{g}}$ of $\cH^k_G$ with a full complement
$\l_1\le\cdots\le\l_{B_k}$ of real eigenvalues (by the symmetry of the
conformally invariant quadratic form). Denoting by $\kappa_{i}$ the
corresponding eigenvectors, the $(\l,\kappa)$ solving our problem
must be one of the pairs $(\l_i,\kappa_i)$ from this list.

\vspace{0.2cm}

\noindent Thomas Branson, Department of Mathematics, The 
University of Iowa, Iowa City IA
52242 USA\newline
{\tt thomas-branson@uiowa.edu}
\vspace{0.2cm}

\noindent A.\ Rod Gover, 
Department of Mathematics,
The University of Auckland,
Private Bag 92019,
Auckland 1,
New Zealand\newline
{\tt gover@math.auckland.ac.nz}

\end{document}